\documentclass[a4paper 11pt]{article}
\usepackage{latexsym,amsmath, amsthm, amssymb}
\catcode`@=11 \long\def\@makefntext#1{\noindent #1}
\newskip\tabcentering \tabcentering=1000pt plus 1000pt minus 1000pt
\def\MCH#1#2{\setbox0=\hbox{\raise#1\hbox{#2}}\smash{\box0}}% move char

\def\@evenfoot{}\def\@oddfoot{}

\def\@evenhead{\hbox to\textwidth{\footnotesize\rm\thepage \hfill
{\it }}}

\def\@oddhead{\hbox to \textwidth{\footnotesize{\it
} \hfill\thepage}}

\def\proof{\vspace{2mm}\noindent{\it Proof}\quad}

\newtheorem{lem}{Lemma}[section]
\newtheorem{thm}{Theorem}[section]
\newtheorem{definition}{Dfinition}[section]

\newtheorem{remark}{Remark}[section]

\newtheorem{pro}{Proposition}[section]

\def\bc{\begin{center}}
\def\ec{\end{center}}

\def\hang{\hangindent\parindent}
\def\textindent#1{\indent\llap{\qquad #1\ \ \enspace}\ignorespaces}
\def\ref{\par\hang\textindent}

\begin{document}
 \abovedisplayskip=8pt plus 1pt minus 1pt
\belowdisplayskip=8pt plus 1pt minus 1pt
%-------------------  First Head  -----------------------------------------
\thispagestyle{empty} \vspace*{-3.0truecm} \noindent
\parbox[t]{6truecm}{\footnotesize\baselineskip=11pt\noindent  {} %Acta Mathematica
%%Sinica, English Series\\
%%1999, Jan., Vol.15, No.1, p. 1--11\\
%%Http://www.ActaMath.com\\
%DOI:
 } \hfill
%%\parbox[c]{6truecm}{\vbox{\hsize 3.6576 true cm %
%%  \vskip 3.8 true cm %1.8373
%%  \relax\hbox to0.4\hsize{\hbox to0pt{\special{BMF=actmark.BMF}}\hss}\hss}}
%\hbox to\textwidth{\vbox{\footnotesize\baselineskip=11pt\noindent
%Acta Mathematica Sinica, English Series\hfill LOGO\\
%1999, Jan., Vol.15, No.1, pp. 1--11\hfill \copyright Spring-Verlag 1999}
%\vbox{\hsize3.6576 true cm
%  \vskip1.8373 true cm
%  \relax\hbox to\hsize{\hbox to0pt{\special{BMF=ACTMARK.BMF}}\hss}\hss}}
%===================Text=============================================
\vspace{1 true cm}

\bc{\large\bf Optimal stopping under $g_\Gamma$-expectation }
%\footnotetext{\footnotesize Received February 24, 1998, Revised
%September 1, 1998, Accepted September 9, 1996}}
\footnotetext{\footnotesize + Corresponding author.\\
} \ec

\vspace*{0.1 true cm}
\bc{\bf Helin Wu$^{+}$  \\
{\small\it School of Mathematics, Shandong University, Jinan 250100, China\\
\small\it \quad E-mail: wuhewlin@gmail.com}}\ec \vspace*{3 true mm}

\begin{abstract}

In this paper, we solve  the  existence problem of optimal stopping
problem under some kind of nonlinear expectation named
$g_\Gamma$-expectation which was recently  introduced in  Peng, S.G.
and Xu, M.Y. [8]. Our method based on our preceding work on the
continuous property of $g_\Gamma$-solution. Generally, the strict
comparison theorem does not hold under such nonlinear expectations
any more, but  we can still modify the classical method to find out
an optimal stopping time via continuous property. The mainly used
theory in our paper is the monotonic limit theorem of BSDE and
nonlinear decomposition theorem of Doob-Meyer's type developed by
Peng S.G. [6]. With help of these useful theories,  a RCLL
modification of the value process can also be obtained by a new
approach instead of down-crossing inequality.

\end{abstract}

\section {Introduction}
Optimal stopping is a very classical and meaningful problem in pure
stochastic analysis and applications, a well-known example is that
the price of American claims under complete market without arbitrage
can be represented by the value function of an optimal stopping
under some linear expectation induced by a probability. In recent
years, nonlinear expectation become more and more wildly studied by
authors. Among all kinds of nonlinear expectations, $g$-expectation
which was introduced by Peng S.G. [7] is a nice example for it
enjoys many nice properties like linear expectation such as
continuous property an strict comparison property as well as
time-consistence property. In more general case, for a given family
of $\sigma$-fields $\{\mathcal{F}_t\}_{t\in[0,T]}$, some kind of
$\mathcal{F}$-expectation can be defined by  axioms.  An interesting
problem that when $\mathcal{F}$-expectation can be represented by
$g$-expectation was considered in Coquet, F., Hu, Y.,  M∩emin, J.,
and Peng, S.G. [2] where  we refer to the definition of
$\mathcal{F}$-expectation.

From an application point of view, many nonlinear expectations are
inevitable just because of the world in our reality is not idea and
perfect. For example, the pricing and hedging problem can be modeled
by linear BSDE under complete market without arbitrage while BSDE
driven by nonlinear generator function $g$ becomes reasonable when
the market is incomplete or other kinds of constraints be
considered.

In the framework of nonlinear expectation.  Riedel, F. [9] studied
the optimal stopping problem with multiple priors. The author
developed a theory of the optimal stopping along the classical lines
using and extending suitable results from martingale theory in
finite discrete time model. This approach works as long as the set
of priors is time consistent.  Kr\"{a}tschmer, V. and  Schoenmakers,
J. [4] considered the optimal stopping for more general dynamic
utility functionals satisfying nice properties such as time
consistency an recursiveness but without strict comparison property
in finite discrete time case. In their paper, the authors provided
sufficient conditions for Bellman principle and the existence of
optimal stopping.
    For continuous  time, an optimal stopping problem was considered
    under ambiguity by  Cheng, X. and Riedel, F. [3]. In that paper, the author
    solve the optimal stopping problem using nonlinear
    Doob-Meyer-Peng decomposition of $g$-supermartingale. However,
    Bayraktar, E. and  Yao, S [1] developed a theory for solving
    continuous time stopping problems for general non-linear
    expectations. Given a stable family of
    $\mathcal{F}$-expectations $\{\mathcal{E}_i\}_{i\in I}$ defined
    well on a common domain, the authors considered the optimal
    problems $$\sup_{(i,\tau)\in I\times
    S_{0,T}}\mathcal{E}_i(Y_\tau+H^i_\tau)\eqno(1.1)$$
 and
    $$\sup_{\tau \in S_{0,T}}\inf_{i\in I}\mathcal{E}_i(Y_\tau+H^i_\tau).\eqno(1.2)$$
where $S_{0,T}$ denotes the whole stopping times valued on $[0,T]$
and $(Y_t + H_t^i),i\in I$ are the model-dependent reward processes.
    Among all above papers except for  Kr\"{a}tschmer, V. and  Schoenmakers, J. [4]
    , nonlinear expectations all satisfy the uncommon property of
    strict comparison and stable property.

    In our paper, we consider the optimal stopping problem under
    $g_\Gamma$-expectation which was introduced by Peng, S.G. and Xu, M.Y.
    [8] as follows:
    $$\sup_{\tau \in S_{0,T}}\mathcal{E}_0^{g,\phi}(X_\tau ).\eqno(1.3)$$
    where $\mathcal{E}_0^{g,\phi}(\cdot )$ is the
    $g_\Gamma$-expectation  under some constraints $\phi(t,y,z)=0$
    well defined on some suitable space
    and $(X_t)$ is a reward process satisfying some mild assumptions.

    Although the  $g_\Gamma$-expectation can not easily be represented  by a
    stable class of $g$-expectations, but it is still a increasing limit of a
    sequence of $g_n$-expectations via penalization method.
    A main difficulty is that the strict comparison theorem may not
    holds for $g_\Gamma$-expectation any more. The method used to
    solve optimal stopping problem in above mentioned papers must be
    modified to work well in our case. Fortunately, with the help of
    some results  about the continuous property of $g_\Gamma$-expectation
   obtained in our preceding paper Wu, H.L. [10], we can still find out an
   optimal solution of this problem.

   Our paper is organized  as follows: In section 2, we give some
   necessary definitions such as $g_\Gamma$-expectation  and some
   useful properties of it. In section 3, we work out the optimal
   problem by a modified method of classical one.

\section {CBSDE  and $g_\Gamma$-expectation }

\noindent Given a probability space $(\Omega,\mathcal{F},P)$ and
$R^d$-valued Brownian motion \mbox{$W(t)$}, we consider a sequence
$\{(\mathcal{F}_t);t\in[0,T]\}$ of filtrations generated by Brownian
motion $W(t)$  and augmented by P-null sets. $\mathcal{P}$  is the
$\sigma$-field of predictable sets of  $\Omega\times[0,T]$. We use
$L^2(\mathcal{F}_T)$ to denote the space of all $F_T$-measurable
random variables $\xi:\Omega\rightarrow R^d$ for which
$$\parallel \xi\parallel^2=E[|\xi|^2]<+\infty.$$
and use $H_T^2(R^d)$ to denote  the space of predictable process
$\varphi:\Omega\times[0,T]\rightarrow R^d$ for which
$$\parallel \varphi\parallel^2=E[\int_0^T|\varphi|^2]<+\infty.$$
For a given probability $P$, we denote the Banach  space of all
P-essentially   bounded real functions  on a probability space
$(\Omega,\mathcal{F_T},P)$ as $L^\infty(\mathcal{F}_T)$.

Given a function $\varphi:[0,T]\times R\times R^d\rightarrow R$,
following assumptions always used in theory of BSDE.
$$|\varphi(\omega,t,y_1,z_1)-\varphi(\omega,t,y_2,z_2)|\leq M(|y_1-y_2|+|z_1-z_2|), \ \ \forall
(y_1,z_1),(y_2,z_2)\eqno(A1)$$ for some $M>0$.

$$\varphi(\cdot,y,z)\in H_T^2(R)\quad \forall y\in R,\, z\in R^d  \eqno(A2)$$

The backward stochastic differential equation (shortly BSDE ) driven
by $g(t,y,z)$ is given by
$$-dy_t=g(t,y_t,z_t)dt-z^*_tdW(t) \eqno(2.1)$$
 where $y_t\in R$ and $W(t)\in R^d$.
Suppose that $\xi\in L^2(\mathcal{F}_T)$ and $g$ satisfies (A1) and
(A2), Pardoux, E., Peng, S.G. [5] proved the existence  of adapted
solution $(y(t),z(t))$ of such BSDE. We call $(g,\xi)$ standard
parameters for the BSDE.

We call the pair $(y_t, z_t)$  satisfying $(2.1)$ a g-solution, but
when an increasing process is added in a BSDE,  the notation of
super-solution is introduced by researchers.

\begin{definition}
(super-solution) A super-solution of a BSDE associated with the
standard parameters $(g,\xi)$ is a vector process $(y_t,z_t,C_t)$
satisfying
$$-dy_t=g(t,y_t,z_t)dt+dC_t-z^*_tdW(t),\quad y_T=\xi,\eqno(2.2)$$
or being equivalent to
$$y_t=\xi+\int_t^Tg(s,y_s,z_s)ds-\int_t^Tz^*_sdW_s+\int_t^TdC_s, \eqno(2.2')$$
where $(C_t,t\in[0,T])$ is an increasing, adapted, right-continuous
process with $C_0=0$ and $z_t^*$ is the transpose  of $z_t$.

\end{definition}

In many analysis and applications, constraints always put on $(y_t,
z_t)$. We formulate the constraints like stated in Peng, S.G. [6].
For a given function $\phi(t,y,z): [0,T]\times R\times
R^d\rightarrow R^+$, we define a subset in $[0,T]\times R\times R^d$
as $\Gamma_t\triangleq\{(t,y,z)|\phi(t,y,z)=0\}$.

 A super-solution $(y_t,z_t,A_t)$ is said to satisfies  constraints if the following
 condition holds,
 $$(t,y_t,z_t)\in \Gamma_t.\eqno(2.3)$$

 \noindent Constraints
like (2.3) is always  considered  in this paper. In such case, we
give the following definition.
\begin{definition}( $g_\Gamma$-solution or the minimal solution ) A g-supersolution $(y_t, z_t, C_t)$ is said to be the
 the minimal solution of a constrained backward differential stochastic equation (shortly CBSDE), given $y_T=\xi$,
 subjected to the constraint $(2.3)$ if for any other g-supersolution
$(y'_t, z'_t, C'_t)$ satisfying $(2.3)$  with $y'_T=\xi$, we have
$y_t\leq y'_t $ a.e., a.s.. The  minimal solution is denoted by
$\mathcal{E}_t^{g,\phi}(\xi)$ and for convenience called as
$g_\Gamma$-solution. Sometimes, we also  call $g_\Gamma$-expectation
$\mathcal{E}_t^{g,\phi}(\xi)\triangleq y_t$ the dynamic
$g_\Gamma$-expectation with constraints $(2.3)$.
\end{definition}
For any $\xi\in L^2(\mathcal{F}_T)$, we denote
$\mathcal{H}^{\phi}(\xi)$ as the set of g-supersolutions
$(y_t,z_t,C_t)$ subjecting to $(2.3)$ with $y_T=\xi$.  When
$\mathcal{H}^{\phi}(\xi)$ is not empty, Peng, S.G. [6] proved that
$g_\Gamma$-solution exists.

In general case, unlike $g$-solution, the increasing part of
$g_\Gamma$-solution is different  with different terminal value and
it is impossible to get a similar priori estimation. The continuous
property seems hard to hold, however,  we can prove it is still
continuous from below similarly like Wu, H.L. [10].
\begin{pro}
Suppose  the generator function $g(t,y,z)$ and the  constraint
function $\phi(t,y,z)$ both satisfy conditions (A1) and (A2),
  $\{\xi_n\in L^\infty(\mathcal{F}_T), \ n=1,2,\cdots\}$ is an norm-bounded increasing
sequence in $L^\infty(\mathcal{F}_T)$  and converges almost surely
to $\xi\in L^\infty(\mathcal{F}_T)$, if
$\mathcal{E}_t^{g,\phi}(\zeta)$ exists for
$\zeta=\xi,\xi_n,n=1,2,\cdots$, then
$$\lim_{n\rightarrow \infty}\mathcal{E}_t^{g,\phi}(\xi_n)= \mathcal{E}_t^{g,\phi}(\xi).$$
\end{pro}

In order to obtain a whole continuity,  we  always assume that both
$g$ and $\phi$ are convex functions.

The convexity of $\mathcal{E}_t^{g,\phi}(\xi)$ can be easily deduced
from the same proposition of solution of BSDE with convex generator
function, see also Peng, S.G. and Xu, M.Y. [8].
\begin{pro}
Suppose that $\phi(t,y,z)$ and $g(t,y,z)$ are both convex in $(y,z)$
and satisfy (A1) and (A2), then
$$\mathcal{E}_t^{g,\phi}(a\xi+(1-a)\eta)\leq a\mathcal{E}_t^{g,\phi}(\xi)+(1-a)\mathcal{E}_t^{g,\phi}(\eta)\quad \forall t\in[0,T]$$
holds for any $\xi,\eta$ in the effective domain of CBSDE and
$a\in[0,1]$.
\end{pro}

\proof According to Peng, S.G. [6], the solutions  $y_t^n(\xi)$ of
$$y_t^n(\xi)=\xi + \int_t^Tg(y_s^n(\xi),z_s^n,s)ds+A_T^n-A_t^n-\int_t^Tz_s^ndW_s.$$
is an increasing sequence and converges  to
$\mathcal{E}_t^{g,\phi}(\xi)$, where $$A_t^n: =
n\int_0^t\phi(y_s^n,z_s^n,s)ds.$$ For any fixed $n$, by the
convexity of $g$ and $\phi$, $y_t^n(\xi)$ is a convex in $\xi$, that
is
$$y_t^n(a\xi+(1-a)\eta)\leq ay_t^n(\xi)+(1-a)y_t^n(\eta),$$
taking limit as $n\rightarrow\infty$, we get the required result.
\hspace*{\fill}$\Box$

By the same method of penalization, we can get the comparison
theorem of $\mathcal{E}_t^{g,\phi}(\xi)$ .
\begin{pro}
Under the same assumptions as above proposition, we have
$$\mathcal{E}_t^{g,\phi}(\xi )\leq
\mathcal{E}_t^{g,\phi}(\eta)$$ for any $\xi, \eta\in L_T^2(R)$ when
$P(\eta\geq\xi)=1$.

\end{pro}

In order to make the domain of definition of CBSDE  more explicitly
for our use, we give another mild assumption below,

$$\varphi(\cdot,y,0)=0 \qquad \forall y\in R.\eqno (A3)$$
The following result can be easily obtained with the help of Peng,
S.G. and Xu, M.Y. [8].

\begin{pro}
Suppose the generator function $g$ and the constraint function
$\phi$ satisfy assumptions $A(i),i=1,2,3$, then the
$g_\Gamma$-solution exists for any $\xi\in L^\infty(\mathcal{F}_T)$
with terminal condition $y_T=\xi$.

\end{pro}
\proof In the paper Peng, S.G. and Xu, M.Y. [8], the author define a
new subspace $L^2_T(R)$:
$$L^2_{+,\infty}(\mathcal{F}_T)\triangleq\{\xi\in
L^2(\mathcal{F}_T), \xi^+\in L^\infty(\mathcal{F}_T)\}.$$

For any $\xi\in L^2_{+,\infty}(\mathcal{F}_T)$ with terminal
condition $y_T=\xi$,  the existence of $g_\Gamma$-solution   was
proved in that paper under the assumption
$$g(t,y,0)\leq L_0 +M|y| \quad \text{and}\quad (y,0)\in \Gamma_t\eqno(2.4)$$
holds for a large constant $L_0$ and for any $ y\geq L_0$,

It is obvious $L^\infty(\mathcal{F}_T)\subset
L^2_{+,\infty}(\mathcal{F}_T)$ and under assumptions $A(i),i=1,2,3$,
$(2.4)$ holds for any $L_0\geq 0$ and $M$ in $(A1)$, thus
$g_\Gamma$-solution is defined well on the whole space
$L^\infty(\mathcal{F}_T)$.\hspace*{\fill}$\Box$.

 The following nice
properties of $\mathcal{E}_t^{g,\phi}(\cdot)$ will be helpful in our
study, their proofs can be found in Peng, S.G. and Xu, M.Y. [8].
\begin{pro}
Suppose the generator function $g$ and the constraint function
$\phi$ satisfy assumptions $A(i),i=1,2,3$, then the
$g_\Gamma$-expectation satisfies:
\begin{itemize}

\item[(i)] Self-preserving: $\mathcal{E}_t^{g,\phi}(\xi_t)=\xi_t$ for
any $\xi_t\in L^\infty(\mathcal{F}_t).$

\item[(ii)] Time
consistency:$\mathcal{E}_s^{g,\phi}(\mathcal{E}_t^{g,\phi}(\xi))=\mathcal{E}_s^{g,\phi}(\xi),\qquad
0\leq s\leq t\leq T \quad \xi\in \in L^\infty(\mathcal{F}_T).$
\item[(ii)] 1-0 law:
$1_A\mathcal{E}_t^{g,\phi}(\xi)=\mathcal{E}_t^{g,\phi}(1_A\xi),\quad
\forall A\in\mathcal{F}_t.$
\end{itemize}

\end{pro}

When both $g(t,y,z)$ and $\phi(t,y,z)$ are convex in $(y,z)$, with
the help of convex analysis, Wu, H.L. [10] has proved
$g_\Gamma$-solution is continuous according to the norm of
$L^\infty(\mathcal{F}_T)$ on the whole space
$L^\infty(\mathcal{F}_T)$. The continuous property will play a
crucial role in our analysis in next section.
\section {Optimal stopping under $g_\Gamma$-expectation }
In this section, we want to find an optimal stopping time which
attains the supermum:
$$\sup_{\tau \in S_{0,T}}\mathcal{E}_0^{g,\phi}(X_\tau ).\eqno(3.1)$$
For simplicity and  making  $\mathcal{E}_0^{g,\phi}(X_\tau )$
meaningful, we assume that the model-depend reward process $X_t,t\in
[0,T]$ is  an adapted, nonnegative process with continuous sample
paths. Furthermore, we still assume  $X_t,t\in [0,T]$  is
 bounded in $L^\infty(\mathcal{F}_T)$. Similarly, like definition in
 Bayraktar, E. and  Yao, S. [1], a process $X_t,t\in [0,T]$ is
 called uniformly-left-continuous if for any sequence
 $\{\tau_n\}_{n\in \mathcal{N}}\subset S_{0,T}$ increasing $a.s $ to $\tau$,
 we can find a subsequence $\{n_k\}_{k\in\mathcal{N}}$ of the set of
 positive nature numbers  $\mathcal{N}$ such that the sequence of
 random variables $X_{\tau_{n_k}},
 k=1,2,\cdots$ converges to $X_\tau$  in $L^\infty(\mathcal{F}_T)$ according to
 norm.

 Under $g_\Gamma$-expectation, we define the value function of the
 optimal stopping problem as
 $$V_t\triangleq  \textrm{ess}\sup_{\tau \in S_{t,T}}\mathcal{E}_t^{g,\phi}(X_\tau
 ).\eqno(3.2)$$
As usual, we define supermartingale (respectively submartingale ,
martingale) under $g_\Gamma$-expectation as done in Peng, S.G. and
Xu, M.Y. [8].
\begin{definition}
A process $(X_t)$ which is adapted  and
$X_t\in\L^\infty(\mathcal{F}_t)$ for every $t\in [0,T]$ is called a
$g_\Gamma$-supermartingale (respectively submartingale , martingale)
on $[0,T]$, if for
  $0\leq s\leq t\leq t\leq T$
 we have
$$\mathcal{E}_s^{g,\phi}(X_t)\leq X_s, (resp,
\geq,=X_s).$$
\end{definition}
Just as classical case, we show that $(V_t)$ defined by $(3.2)$ is a
$g_\Gamma$-supermartingale, it is based on the continuous property
of $g_\Gamma$-solution and the following lemma.

\begin{lem}
For all $t\geq 0$, the family
$$\{\mathcal{E}_t^{g,\phi}(X_\tau):\tau\geq t\}$$
is upwards directed.
\end{lem}
\proof Thanks to the useful property of 1-0 law of
$g_\Gamma$-expectation, we can prove this result by the usual way,
for details, see for example Lemma B.1 in Cheng, X. and Riedel, F.
[3]. \hspace*{\fill}$\Box$

With the help of this lemma, we have
\begin{pro}
Under the assumptions on the reward process in our paper, the value
function $(V_t)$ defined by $(3.2)$ is a $g_\Gamma$-supermartingale.
\end{pro}
\proof For every $t\geq 0$, the lemma above allows us to choose a
sequence $\{\tau_n(t),n=1,2,\cdot\}$ of stopping times greater or
equal $t$ with
$$\mathcal{E}_t^{g,\phi}(X_{\tau_n(t)})\uparrow V_t.$$
Since $\mathcal{E}_t^{g,\phi}(X_{\tau_n(t)})$  converges to
  $V_t$ increasingly, by the continuous property from below of
proposition $(2.1)$ and time consistency  property $(ii)$ in
proposition 2.5,  for $0\leq s \leq t \leq T$, we have
\begin{eqnarray*}
\mathcal{E}_s^{g,\phi}(V_t) &=&\mathcal{E}_s^{g,\phi}(\sup_{\tau \in
S_{t,T}}\mathcal{E}_t^{g,\phi}(X_\tau))=\mathcal{E}_s^{g,\phi}(\lim_{n\rightarrow
\infty} \mathcal{E}_t^{g,\phi}(X_{\tau_n(t)})) \\
&=&\lim_{n\rightarrow
\infty}\mathcal{E}_s^{g,\phi}(\mathcal{E}_t^{g,\phi}(X_{\tau_n(t)}))
=\lim_{n\rightarrow \infty}\mathcal{E}_s^{g,\phi}(X_{\tau_n(t)})\\
&\leq&\textrm{ess}\sup_{\tau \in
S_{s,T}}\mathcal{E}_s^{g,\phi}(X_\tau
 )=V_s.
\end{eqnarray*}
  \hspace*{\fill}$\Box$

To obtain an optimal stopping time, we want to show that there is a
right-continuous modification of $(V_t)$. However, this time, the
strict comparison theorem does not hold for $G_\Gamma$-expectation
anymore in general, so the usual way to find a right-continuous
modification of the value function $(V_t)$ by downcrossing
inequality may not work. Fortunately, with the help of important
results obtained in Peng, S.G. [6], we can still have the following
claim.
\begin{thm}
Under the assumptions in our paper, the value process $(V_t)$
defined by $(3.2)$ has a right-continuous modification.
\end{thm}
\proof Let $g_n=g+n\phi$ as in proposition $(2.2)$, we define the
value function $(V_n(t))$ under $g_n$-expectation
$$V_n(t)\triangleq  \textrm{ess}\sup_{\tau \in S_{t,T}}\mathcal{E}_t^{g_n}(X_\tau
 ).\eqno(3.3)$$
According to  Lemma F.1 in Cheng, X. and Riedel, F. [3] or Lemma 5.2
in Coquet, F., Hu, Y.,  M∩emin, J., and Peng, S.G. [2], $(V_n(t))$
is a $g_n$-supermartingale with a right-continuous modification for
any $n$. By the comparison them of BSDE, without lose of generality,
we can say $(V_n(t))$ is also a RCLL $g$-supermartingale, hence by
Theorem 3.3 in Peng, S.G. [6], it is a $g$-supersolution. At the
same time, we can easily prove that
$$V_n(t)\uparrow V_t. $$

In fact, since $(g_n)$ is an increasing sequence of generator
functions, we have

\begin{eqnarray*}
V_t &=&\textrm{ess}\sup_{\tau \in
S_{t,T}}\mathcal{E}_t^{g,\phi}(X_\tau
 )= \textrm{ess}\sup_{\tau \in S_{t,T}}\sup_{n\in
 \mathcal{N}}\mathcal{E}_t^{g_n}(X_\tau)\\
 &=& \sup_{n\in \mathcal{N}}\textrm{ess}\sup_{\tau \in
 S_{t,T}}\mathcal{E}_t^{g_n}(X_\tau
 )=\sup_{n\in \mathcal{N}}V_n(t).
\end{eqnarray*}
All the process in our paper are  bounded in
$L^\infty(\mathcal{F}_T)$, with the help of Theorem 3.6 of Peng,
S.G. [6], $(V_t)$ is also a RCLL $g$-supersolution or
$g$-supermartingale. \hspace*{\fill}$\Box$

As usual, it is easy to  see that $(V_t)$ is the smallest
$g_\Gamma$-supermartingale with RCLL sample path which we state it
 as a proposition below.
\begin{pro}
The value function process  $(V_t)$ is the smallest
$g_\Gamma$-supermartingale with RCLL sample path which dominates the
reward process $X_t$.
\end{pro}
\proof Suppose $(S_t)$ is another RCLL $g_\Gamma$-supermartingale
with $S_t\geq X_t$ for all $t\in [0,T]$ and $S_\tau$ is meaningful
for any stopping time.

Choose a sequence of stopping times $\{\tau_n(t),n=1,2,\cdot\}$ in
$S_{t,T}$ as in the proof of proposition 3.1, then we have
$$V_t=\lim \mathcal{E}_t^{g,\phi}(X_{\tau_n(t)})\leq \liminf
\mathcal{E}_t^{g,\phi}(S_{\tau_n(t)})\leq S_t.$$

 Hence, $(V_t)$ is the smallest RCLL $g_\Gamma$-supermartingale
 dominating $X$.\hspace*{\fill}$\Box$

 With these results in hand, we then go on to find an optimal
stopping of problem $(3.1)$ by a similar constructive way as usual.

For any $0<\lambda <1, 0\leq t \leq T$, we define the stopping times
$$\tau^\lambda(t) \triangleq \inf \{u\geq t|X_u\geq \lambda V_u\}.$$

The next lemma is a crucial step to construct an optimal stopping
time for our problem.
\begin{lem} With the notation introduced above,
  then we have
$$V_t=\mathcal{E}_t^{g,\phi}(V_{\tau^\lambda(t)}).$$

\end{lem}

\proof Introduce the process
$$W_t\triangleq\mathcal{E}_t^{g,\phi}(V_{\tau^\lambda(t)}).$$
and for each $n$,
$$W_n(t)\triangleq\mathcal{E}_t^{g_n}(V_{\tau^\lambda(t)}),n=1,2,\cdots.$$

Since $(V_t)$ is $g_n$-supermartingale for each $n=1,2,\cdots$, by
the same way similar with Cheng, X. and Riedel, F. [3], $W_n$ is a
$g_n$-supermartingale with RCLL sample paths. Furthermore, we can
claim that $W$ is a  $g_\Gamma$-supermatingale. For $0\leq t\leq
t+u\leq T$ we have
$$\mathcal{E}_t^{g,\phi}(W_{t+u})=\mathcal{E}_t^{g,\phi}(\mathcal{E}_{t+u}^{g,\phi}(V_{\tau^\lambda(t+u)}))
=\mathcal{E}_t^{g,\phi}(V_{\tau^\lambda(t+u)})\leq
\mathcal{E}_t^{g,\phi}(V_{\tau^\lambda(t)})=W_t.
$$

It is obviously that $W_n(t)$ converges increasingly to $W_t$ and
they are all $g$-supermartingales, hence we can use the same skill
in the proof of Theorem (3.2) to prove that $W$ admits a RCLL
modification.

The following proof can go on similarly as the last part of the
proof of Lemma B.3 in Cheng, X. and Riedel, F. [3]. For convenience,
we state it still in our paper.

Let
$$Y_t=\lambda V_t+(1-\lambda)W_t.$$

We claim that $W$ dominates $X$. For $X_t\geq V_t$, we have
$\tau^\lambda(t)=t$, hence $W_t=V_t$ and  $Y_t=V_t\geq X_t$. If
$X_t<\lambda V_t$, we have $W_t\geq 0$ as $X_t\geq 0$, so
$$Y_t\geq\lambda V_t\geq X_t.$$
From Proposition 3.2, we get $Y\geq V$. This equivalent to $W\geq
V$. On the other hand, by definition of $W$ and
$g_\Gamma$-supermartingale property of $V$: $W\leq V$. So we
conclude $W=V$. In other words, we finally get
$$V_t=W_t=\mathcal{E}_t^{g,\phi}(V_{\tau^\lambda(t)}).$$
\hspace*{\fill}$\Box$

Now let us back to the definition of stopping times of
$\tau^\lambda(t)$ at $t=0$, which, for simplicity, we denote it as
$\tau^\lambda$.

Noting that $\tau^\lambda$ is increasing with $\lambda$ and
dominated by the stopping time $\tau^*\triangleq \inf\{t\geq
0:X_t=V_t\}$, we can choose a sequence of real numbers
$(\lambda_n)\subset (0,1)$ such that $\tau^{\lambda_n}$ converges
increasingly to some stopping time $\bar{\tau}$.

We state our last result in this paper.
\begin{thm}
Under our assumptions in our paper about the generator function $g$
and constraint function $\phi$ as well as the model-dependent reward
process $(X_t)$, if furthermore both $g$ and $\phi$ are convex, then
with notations above, the stopping time $\bar{\tau}$ is an optimal
stopping for problem (3.1).
\end{thm}
\proof First, by Lemma 3.2, with $t=0$, we have
$$V_0=\mathcal{E}_0^{g,\phi}(V_{\tau^{\lambda_n}})\leq
\frac{1}{\lambda_n}\mathcal{E}_0^{g,\phi}(X_{\tau^{\lambda_n}}).$$

On the other hand, since $\mathcal{E}_0^{g,\phi}(\cdot)$ is a convex
functional defined well on $L^\infty(\mathcal{F}_T)$, it is then
continuous on the space $L^\infty(\mathcal{F}_T)$ according to the
norm, for details see Wu, H.L. [10].

Our assumption help us to obtain a subsequence of
$(X_{\tau^{\lambda_n}})$, which we still denote as
$(X_{\tau^{\lambda_n}})$, converges to $X_{\bar{\tau}}$ under norm
of $L^\infty(\mathcal{F}_T)$, thus
$$V_0\leq
\lim_{\lambda_n\rightarrow
1}\frac{1}{\lambda_n}\mathcal{E}_0^{g,\phi}(X_{\tau^{\lambda_n}})
=\mathcal{E}_0^{g,\phi}(X_{\bar{\tau}})\leq V_0.$$

Thus $\bar{\tau}$ is an optimal stopping time. \hspace*{\fill}$\Box$

\begin{remark}

 According to Peng, S.G. and Xu, M.Y. [8], $V_t$ can be viewed as
 the solution  of the Reflected Backward stochastic differential
 equation with constraint, and $\tau^*$ defined above is another
 optimal stopping time with the stopped process $(V_{t\wedge
 \tau^*})$ be a $g_\Gamma$-martingale. However, different from
 classical case, $\bar{\tau}$ may not coincide with $\tau^*$, so
 whether the stopped process $(V_{t\wedge \bar{\tau}})$ is also a
 $g_ \Gamma$-martingale is questioned.
\end{remark}
\begin{remark}
By the penalization method to obtain the $g_\Gamma$-solution, $V_t$
can be represented by
$$V_t=\textrm{ess}\sup_{\tau \in
S_{t,T}}\sup_{n}\mathcal{E}_t^{g_n}(X_\tau)$$ which is a stopper and
controller problem.
\end{remark}

\end{document}